\documentclass[leqno,12pt]{article}  
\usepackage{amsmath}
\usepackage{amssymb}

\usepackage[german,english]{babel}
\usepackage{amsthm}

\title{On $\psi$-lattices in modular $(\varphi,\Gamma)$-modules}

\author{\textsc{Elmar Grosse-Kl\"onne}}
\date{}

\theoremstyle{plain} 
\newtheorem{satz}{Theorem}[section]  
\newtheorem{lem}[satz]{Lemma}  
\newtheorem{pro}[satz]{Proposition}  
 
\theoremstyle{remark}

\theoremstyle{definition}

\newcommand{\0}{\ensuremath{\overrightarrow{0}}}

\begin{document}

\maketitle

%


\begin{abstract} Let $F/{\mathbb Q}_p$ be a finite field extension, let $k$ be a finite field extension of the residue field of $F$. Generalizing the $\psi$-lattices which Colmez constructed in \'{e}tale $(\varphi,\Gamma)$-modules over $k[[t]][t^{-1}]$, we define, study and exemplify $\psi$-lattices in \'{e}tale $(\varphi,\Gamma)$-modules over $k[[t_1,\ldots,t_d]][\prod_it_i^{-1}]$ for arbitrary $d\in{\mathbb N}$.
\end{abstract}

\begin{center}{\Large{\bf Introduction}}\end{center} 

Let $F/{\mathbb Q}_p$ be a finite field extension. It was a fundamental insight of Fontaine \cite{font} that ($p$-adically continuous) representations of the absolute Galois group ${\rm Gal}(\overline{F}/F)$ of $F$ on finite free modules over ${\mathbb Q}_p$ or ${\mathbb Z}_p$ or ${\mathbb F}_p$ --- for short: $p$-adic Galois representations --- can equivalently be described by linear algebra objects which he called \'{e}tale $(\varphi,\Gamma)$-modules. These are modules over certain Laurent series rings, in one variable, endowed with commuting semilinear actions by a Frobenius operator ${\varphi}$ and the group $\Gamma={\mathcal O}_F^{\times}$.\footnote{In fact, $\Gamma$ was taken to be a slightly different group in \cite{font}, but there is no substantial difference to the point of view taken here.} Approaching Galois representations through their associated \'{e}tale $(\varphi,\Gamma)$-modules has proven to be an extremely powerful method in numerous contexts. We mention here only the important role which it plays in Colmez' work \cite{colhaupt} on the $p$-adic local Langlands program. Also, the theory of $(\varphi,\Gamma)$-modules has been vastly generalized since into numerous directions. Among these generalizations is the work by Z\'{a}br\'{a}di \cite{zabgal} who showed that (for $F={\mathbb Q}_p$), the representations of the $d$-fold self product ${\rm Gal}(\overline{F}/F)\times\cdots\times{\rm Gal}(\overline{F}/F)$ are in category equivalence with \'{e}tale $(\varphi,\Gamma)$-modules over certain Laurent series rings in $d$ variables (multivariable \'{e}tale $(\varphi,\Gamma)$-modules).

In this note we restrict attention to $p$-modular coefficients only, i.e. our Galois representations (which however remain entirely in the background) and \'{e}tale $(\varphi,\Gamma)$-modules are (in particular) ${\mathbb F}_p$-vector spaces. In this context, an \'{e}tale $(\varphi,\Gamma)$-module is a finite dimensional $k((t))$-vector space ${\bf D}$, for a finite extension $k$ of the residue field of $F$, endowed with said actions by ${\varphi}$ and $\Gamma$. A critical ingredient in the aforementioned work of Colmez was the detection and study of finite $k[[t]]$-lattices spanning ${\bf D}$, stable under $\Gamma$ and a certain operator $\psi$ left inverse to $\varphi$, on which $\psi$ in fact acts surjectively. Among these lattices he identified a minimal one, denoted by ${\bf D}^{\natural}$, and a maximal one, denoted by ${\bf D}^{\sharp}$.

The purposes of this note are the following. Firstly, we want to explain that the study of ${\bf D}^{\natural}$ and ${\bf D}^{\sharp}$ makes sense similarly in the context of ($p$-modular) {\it multivariable} \'{e}tale $(\varphi,\Gamma)$-modules. Secondly, we want to advocate an approach towards these lattices ${\bf D}^{\natural}$ and ${\bf D}^{\sharp}$ in which we rather construct and study their $k$-linear duals. This approach goes hand in hand with a method for constructing \'{e}tale $(\varphi,\Gamma)$-modules. In this latter function, i.e. as a tool for the explicit construction and description of \'{e}tale $(\varphi,\Gamma)$-modules (which typically is quite delicate, since e.g. testing if candidates for $\varphi$- and $\Gamma$-actions, given by explicit power series, really commute with each other can be quite challenging), the method was introduced by Colmez in \cite{colhaupt} and then used later in \cite{elmarsusi}, both in the one-variable case. We intend to use the constructions presented here (the construction and analysis of ${\bf D}^{\natural}$ and ${\bf D}^{\sharp}$) to generalize the work \cite{elmarsusi} to a multivariable setting in the future. Thirdly, by discussing several examples we try to shed some more light on the behaviour of ${\bf D}^{\natural}$ and ${\bf D}^{\sharp}$. Actually, most of these examples pertain to the one-variable case, yet we think that they demonstrate some features of ${\bf D}^{\natural}$ and ${\bf D}^{\sharp}$ which at least have not yet been documented in the literature (although undoubtly known to the experts). \\

{\it Acknowledgements:} I am very grateful for the invitation to the superb event dedicated to perfectoid spaces at the ICTS Bangalore in 2019, and for the opportunity to contribute to its aftermath by means of the present note. I thank the anonymous referee for his helpful remarks on the text.

\section{Multivariable modular \'{e}tale $(\varphi_{\bullet},\Gamma_{\bullet})$-modules} 

{\bf Notations:} Let $F/{\mathbb Q}_p$ be a finite field extension. Denote by $q$ the number of
elements of the
residue field ${\mathbb F}_q$ of $F$. Let $\pi$ be a uniformizer in the
ring of integers ${\mathcal O}_F$ of $F$. Let $k$ be a finite extension field of ${\mathbb F}_q$. Put $\Gamma={\mathcal O}_F^{\times}$.

There is a unique (up to isomorphism) Lubin-Tate group for $F$ with respect to $\pi$. Fixing a coordinate $t$ we write $\Phi(t)$ for the corresponding Lubin-Tate formal power series describing multiplication by $\pi$. Equivalently, to any power series $\Phi(t)\in{\mathcal O}_F[[t]]$ with $\Phi(t)\equiv\pi t$ modulo $t^2{\mathcal O}_F[[t]]$ and $\Phi(t)\equiv t^q$ modulo $\pi{\mathcal O}_F[[t]]$ is associated a Lubin-Tate formal group law, with $\Phi(t)$ representing multiplication by $\pi$, and the resulting formal group (with multiplication by ${\mathcal O}_F$) is independent (up to isomorphism) on the specific $\Phi(t)$. For $\gamma\in\Gamma$ let
$[\gamma]_{\Phi}(t)\in{\mathcal O}_F[[t]]$ denote the power series describing
the action of $\gamma$ in the Lubin-Tate group. Let $D$ be a finite set, and for each $d\in D$ let $t_d$ be a free variable. Put $${k}[[t_\bullet]]={{k}}[[t_d]]_{d\in
  D},\quad\quad\quad\quad {{k}}((t_{\bullet}))={{k}}[[t_{\bullet}]][t_D^{-1}]\quad\mbox{ with }t_D=\prod_{d\in
  D}t_d.$$For each $d\in D$ let $\Gamma_d$ be a copy of $\Gamma$. For
$\gamma\in\Gamma$ let $\gamma_d$ denote the element in $\Gamma_{\bullet}=\prod_{d\in
  D}\Gamma_d$ whose $d$-component is $\gamma$ and whose other components are trivial. The formulae $$\gamma_{d}(t_{d})=[\gamma]_{\Phi}(t_{d}),\quad\quad\quad\quad \gamma_{d_1}(t_{d_2})=t_{d_2}$$with $\gamma\in\Gamma$ and $d, d_1, d_2\in D$ such that
$d_1\ne d_2$ define an action of $\Gamma_{\bullet}$ on ${{k}}((t_{\bullet}))$. Consider the ${k}[[t_\bullet]]$-algebra$${{k}}[[t_{\bullet}]][\varphi_{\bullet},\Gamma_{\bullet}]={{k}}[[t_{d}]]_{d\in D}[\varphi_{d},\Gamma_{d}]_{d\in D}$$ with
commutation rules given by$$x_{d_1}\cdot y_{d_2}=y_{d_2}\cdot
x_{d_1},$$$$\gamma_{d}\cdot
\varphi_{d}=\varphi_{d}\cdot\gamma_{d},\quad\quad \gamma_{d}\cdot
t_{d}=(\gamma_{d}(t_{d}))\cdot\gamma_{d},\quad\quad\varphi_{d}\cdot
t_{d}=\Phi(t_{d})\cdot \varphi_{d}=t_d^q\cdot\varphi_d$$for $\gamma\in\Gamma$ and $x,y\in\Gamma\cup\{\varphi,t\}$
and $d, d_1, d_2\in D$ with $d_1\ne d_2$. Similarly we define the $k((t_{\bullet}))$-algebra
${{k}}((t_{\bullet}))[\varphi_{\bullet},\Gamma_{\bullet}]$ and its subalgebra ${{k}}((t_{\bullet}))[\Gamma_{\bullet}]$.\\

{\bf Definition:} An \'{e}tale $\varphi_{\bullet}$-module over
${{k}}((t_{\bullet}))$ is a ${{k}}((t_{\bullet}))[\varphi_{\bullet}]$-module ${\bf D}$
which is finitely generated over ${{k}}((t_{\bullet}))$ such that for each $d\in D$ the linearized structure map$${\rm id}\otimes\varphi_{d}:{{k}}((t_{\bullet}))\otimes_{\varphi_{d},{{k}}((t_{\bullet}))} {\bf
  D}\longrightarrow{\bf
  D}$$is bijective. An \'{e}tale $(\varphi_{\bullet},\Gamma_{\bullet})$-module over
${{k}}((t_{\bullet}))$ is a
${{k}}((t_{\bullet}))[\varphi_{\bullet},\Gamma_{\bullet}]$-module whose
underlying $\varphi_{\bullet}$-module is \'{e}tale. In the case $|D|=1$ we
drop the indices $(.)_{d}$ resp. $(.)_{\bullet}$ and simply talk about \'{e}tale $(\varphi,\Gamma)$-modules over
${{k}}((t))$.\\

{\bf Remark:} The action of $\Gamma_{\bullet}$ on an \'{e}tale $(\varphi_{\bullet},\Gamma_{\bullet})$-module is automatically continuous for the weak topology.

\begin{lem}\label{review1} (\cite{egk} Lemma 4, Proposition 6)

  The category of \'{e}tale
  $(\varphi_{\bullet},\Gamma_{\bullet})$-modules over ${{k}}((t_{\bullet}))$ is abelian.

  Regard ${k}((t))$ as a ${{k}}((t_{\bullet}))$-module by means of $t_d\cdot x=tx$ for $d\in D$ and $x\in {k}((t))$. The functor $${\bf
  D}\mapsto {{k}}((t))\otimes_{{{k}}((t_{\bullet}))}{\bf
  D}$$ is
  an exact and faithful functor from the category of \'{e}tale $(\varphi_{\bullet},\Gamma_{\bullet})$-modules over ${{k}}((t_{\bullet}))$ to the category of \'{e}tale $(\varphi,\Gamma)$-modules over ${{k}}((t))$.\footnote{We will not make use of the second statement of Lemma \ref{review1} in the following.}
\end{lem}

\begin{satz}\label{sosego} (a) (Fontaine \cite{font}, Kisin-Ren \cite{kisren}, Schneider \cite{peterlec}) There is an equivalence between the  category of \'{e}tale
  $(\varphi,\Gamma)$-modules over ${{k}}((t))$ and the category of continuous representations of ${\rm Gal}(\overline{F}/F)$ on finite dimensional $k$-vector spaces.

  (b) (Z\'{a}br\'{a}di \cite{zabgal}) Assume $F={\mathbb Q}_p$. There is an equivalence between the  category of \'{e}tale
  $(\varphi_{\bullet},\Gamma_{\bullet})$-modules over ${{k}}((t_{\bullet}))$
  and the category of continuous representations of ${\rm
    Gal}(\overline{{\mathbb Q}_p}/{\mathbb Q}_p)\times\cdots\times {\rm
    Gal}(\overline{{\mathbb Q}_p}/{\mathbb Q}_p)$ (with $|D|$ many factors
  indexed by $D$) on finite dimensional $k$-vector spaces.

(c) Statements (a) and (b) merge into a common generalization (Pupazan \cite{pupazan}): For any $F$, there is an equivalence between the  category of \'{e}tale
  $(\varphi_{\bullet},\Gamma_{\bullet})$-modules over ${{k}}((t_{\bullet}))$
and the category of continuous representations of ${\rm
  Gal}(\overline{F}/F)\times\cdots\times {\rm
  Gal}(\overline{F}/F)$ on finite dimensional $k$-vector spaces.
  
\end{satz}

{\bf Definition:} A $\psi$-operator on ${{k}}[[t_{\bullet}]]$ is a system
$\psi_{\bullet}=(\psi_d)_{d\in D}$ of additive maps $$\psi_{d}:{{k}}[[t_{\bullet}]]\longrightarrow
{{k}}[[t_{\bullet}]]$$ for $d\in D$ such that
$\psi_{d_1}(\gamma_{d_2}(t_{d_1}))=\gamma_{d_2}(\psi_{d_1}(t_{d_1}))$
for all $\gamma\in\Gamma$ and $d_1, d_2\in D$, such that $\psi_{d_1}(t_{d_2})=t_{d_2}$ for $d_1\ne
d_2$ and such
that the following holds true: If we view the $\varphi_{d}$ as acting on
${{k}}[[t_{\bullet}]]$, then $\psi_{d_1}\circ\varphi_{d_2}=\varphi_{d_2}\circ\psi_{d_1}$ for $d_1\ne
d_2$, but $$\psi_{d}(\varphi_{d}(a)x)=a\psi_{d}(x)$$ for $a,x\in {{k}}[[t_{\bullet}]]$.\footnote{Notice that we do not require $\psi_{d}(1)=1$.} 

\begin{lem}\label{erleicht} (\cite{egk} Lemma 3) There is a $\psi$-operator $\psi_{\bullet}$ on ${{k}}[[t_{\bullet}]]$ such that each $\psi_{d}$ is surjective.  
\end{lem}

To be explicit, in the case where $|D|=1$ (the general case is handled factor by factor) and  $\Phi(t)=\pi t+t^q$, we may choose $\psi_{k((t))}$ on $k((t))$ such that for $m\in{\mathbb Z}$ and
$0\le i\le q-1$ we have\footnote{Notice that $\frac{q}{\pi}=0$ (in $k$) if $F\ne{\mathbb Q}_p$.}\begin{gather}\psi_{k((t))}(t^{mq+i})=\left\{\begin{array}{l@{\quad:\quad}l}\frac{q}{\pi}t^m& i=0\\0
&  1\le i\le q-2\\t^m& i=q-1\end{array}\right.\label{psieqqp}.\end{gather}

In the following, we fix $\psi_{\bullet}$ as in Lemma \ref{erleicht}.\\

Let ${\bf D}$ be an
\'{e}tale $(\varphi_{\bullet},\Gamma_{\bullet})$-module over
$k((t_{\bullet}))$. For $d\in D$ we define the composed map$$\psi_{d}:{\bf D}\longrightarrow k((t_{\bullet}))\otimes_{\varphi_{d},k((t_{\bullet}))}{\bf D}\longrightarrow{\bf D}$$where the first arrow is the inverse of the structure isomorphism ${\rm id}\otimes\varphi_{d}$, and where the second arrow is given by $a\otimes x\mapsto \psi_{d}(a)x$. \\

\begin{lem} For all $x\in {\bf D}$, $\gamma\in\Gamma$, $a\in k((t_{\bullet}))$ and $d\in D$ we have$$\psi_d(a\varphi_d(x))=\psi_d(a)x,\quad\quad \psi_d(\varphi_d(a)x)=a\psi_d(x),\quad\quad  \gamma_{d}(\psi_d(x))=\psi_d(\gamma_{d}(x)).$$For all $x\in {\bf D}$ and $d_1\ne d_2\in D$ we have $$\psi_{d_1}(\psi_{d_2}(x))=\psi_{d_2}(\psi_{d_1}(x)),\quad\quad \psi_{d_1}(\varphi_{d_2}(x))=\varphi_{d_2}(\psi_{d_1}(x)).$$
\end{lem}

{\sc Proof:} The formula $\psi_d(a\varphi_d(x))=\psi_d(a)x$ is immediate from
the construction. To see the formula $\psi_d(\varphi_d(a)x)=a\psi_d(x)$, write
$x=\sum_ia_i\varphi_d(e_i)$ with $e_i\in {\bf D}$ and $a_i\in
k((t_{\bullet}))$ (this is possible as ${\bf D}$ is \'{e}tale). We then
compute$$\psi_d(\varphi_d(a)x)=\sum_i\psi_d(\varphi_d(a)a_i\varphi_d(e_i))=\sum_i\psi_d(\varphi_d(a)a_i)e_i$$$$=a\sum_i\psi_d(a_i)e_i=a\sum_i\psi_d(a_i\varphi_d(e_i))=a\psi_d(x).$$To
see the formula $\gamma_{d}(\psi_d(x))=\psi_d(\gamma_{d}(x))$ observe that,
since the actions of $\gamma_{d}$ and $\varphi_d$ on $k[[t_{\bullet}]]$ commute, and since
$\Gamma_{\bullet}$ acts semilinearly on ${\bf D}$, the additive map
$$k((t_{\bullet}))\otimes_{\varphi_d,k((t_{\bullet}))}{\bf D}\to k((t_{\bullet}))\otimes_{\varphi_d,k((t_{\bullet}))}{\bf
  D},$$$$a\otimes b\mapsto \gamma_d(a)\otimes\gamma_d(b)$$ is the map corresponding
to the action of $\gamma_d$ on ${\bf D}$ under the isomorphism ${\rm id}\otimes\varphi_d$, and under
$a\otimes x\mapsto \psi_d(a)x$ it commutes with $\gamma_d$ acting on ${\bf D}$ since the actions of $\gamma_d$
and $\psi_d$ on $k((t_{\bullet}))$ commute. The remaining commutation formulae
are clear.\hfill$\Box$\\

{\bf Definition:} For a $k$-vector space $\Delta$ we write $\Delta^*={\rm Hom}_k(\Delta,k)$. We say that a ${{k}}[[t_{\bullet}]]$-module $\Delta$ is admissible if it is a torsion module\footnote{In \cite{egk} we had not included this torsion condition into the definition of admissibility; however, for staying consistent with established terminology we should have done so. (Yet, the omission of this condition in the paper \cite{egk} does not invalidate its results.)} over ${{k}}[[t_{d}]]$ for each $d\in D$ and if $$\Delta[t_{\bullet}]=\{x\in\Delta\,|\,t_dx=0\mbox{ for each }d\in D\}$$is a finite dimensional $k$-vector space. \\

\begin{pro}\label{nopsi} (\cite{egk} Proposition 5) Let $\Delta$ be a finitely
  generated ${{k}}[[t_{\bullet}]][\varphi_{\bullet},\Gamma_{\bullet}]$-module which is admissible as a $k[[t_{\bullet}]]$-module and satisfies $\Delta={{k}}[[t_{\bullet}]]\varphi_d(\Delta)$ for each $d\in D$. Then
  $\Delta^*\otimes_{{{k}}[[t_{\bullet}]]}{{k}}((t_{\bullet}))$ is
  in a natural way an \'{e}tale $(\varphi_{\bullet},\Gamma_{\bullet})$-module over
${{k}}((t_{\bullet}))$. The functor $\Delta\mapsto
  \Delta^*\otimes_{{{k}}[[t_{\bullet}]]}{{k}}((t_{\bullet}))$ is exact.
\end{pro}

We remark that the ${k}[[t_{\bullet}]][\Gamma_{\bullet}]$-action on
$\Delta^*\otimes_{{{k}}[[t_{\bullet}]]}{{k}}((t_{\bullet}))$ results from the
${k}[[t_{\bullet}]][\Gamma_{\bullet}]$-action on $\Delta^*$ given by the formulae$$(a\cdot\ell)(\delta)=\ell(a\delta),$$$$(\gamma\cdot \ell)(\delta)=\ell(\gamma^{-1}\delta)$$for $a\in
{k}[[t_{\bullet}]]$, $\ell\in \Delta^*$, $\delta\in\Delta$ and
$\gamma\in\Gamma_{\bullet}$. The $\varphi_{\bullet}$-action on
$\Delta^*\otimes_{{{k}}[[t_{\bullet}]]}{{k}}((t_{\bullet}))$ is a certain right
inverse to the dual of the $\varphi_{\bullet}$-action on $\Delta$; its
construction involves the $\psi$-operator $\psi_{\bullet}$ on ${{k}}[[t_{\bullet}]]$.\\

\begin{lem}\label{villpasteur} Let $\Delta$ and ${\bf
    D}=\Delta^*\otimes_{k[[t_{\bullet}]]}k((t_{\bullet}))$ be as in Proposition
  \ref{nopsi} and suppose that each $t_{d}$ acts surjectively on $\Delta$.

(a) The natural map $\Delta^*\to{\bf
    D}=\Delta^*\otimes_{k[[t_{\bullet}]]}k((t_{\bullet}))$ is injective. Each
  $\psi_{d}$ respects $\Delta^*$ and acts on it by the rule $$[\ell:\Delta\to
    k]\quad\mapsto\quad [\Delta\to k, x\mapsto\ell(\varphi_{d}(x))].$$

(b) If each
  $\varphi_{d}$ acts injectively on $\Delta$ then each $\psi_{d}$ acts
  surjectively on $\Delta^*$.  
\end{lem} 

{\sc Proof:} This follows immediately from the construction (given in \cite{egk}
Proposition 5) referred to in Proposition
\ref{nopsi}.\hfill$\Box$\\

\section{The lattices  ${\bf D}^{\natural}$ and ${\bf D}^{\sharp}$ inside ${\bf D}$}

We write $\psi_D=\prod_{d\in D}\psi_d$ and $\varphi_D=\prod_{d\in D}\varphi_d$ (as $k$-linear operators on $k((t_{\bullet}))$).

Let ${\bf D}$ be an
\'{e}tale $(\varphi_{\bullet},\Gamma_{\bullet})$-module over
$k((t_{\bullet}))$. We
call a finitely generated $k[[t_{\bullet}]]$-submodule of ${\bf D}$ a lattice
in ${\bf D}$ if it generates ${\bf D}$ (as a $k((t_{\bullet}))$-module). 
  
\begin{lem}\label{codoni} Let $E$ be a lattice in ${\bf D}$, let $d\in D$.

(a) $\psi_d(E)$ is a $k[[t_{\bullet}]]$-module.

(b) If $\varphi_d(E)\subset E$ then $E\subset\psi_d(E)$.

(c) If $E\subset k[[t_{\bullet}]]\cdot \varphi_d(E)$ then $\psi_d(E)\subset E$.

(d) If $\psi_d(E)\subset E$ then $\psi_d(t_d^{-1}E)\subset t_d^{-1}E$. For each $x\in {\bf D}$ there is some $n(x)\in{\mathbb N}$ such that for all $n\ge n(x)$ we have $\psi_D^n(x)\in t_D^{-1}E$.

\end{lem}

{\sc Proof:} (a) Use $\psi_d(\varphi_d(a)x)=a\psi_d(x)$ for $a\in
k((t_{\bullet}))$ and $x\in {\bf
  D}$.

(b) Choose $a\in k[[t_{\bullet}]]$ with $\psi_d(a)=1$. For $e\in E$ we have $e=\psi_d(a\varphi_d(e))$ which belongs to $\psi_d(E)$ since $\varphi_d(E)\subset E$.

(c) Let $e\in E$. By assumption there are $e_i\in E$ and $a_i\in k[[t_{\bullet}]]$ with $e=\sum_ia_i\varphi_d(e_i)$, hence $\psi_d(e)=\sum_i\psi_d(a_i)e_i\in E$.

(d) For $i\ge1$ we
have $$\psi_d(\varphi_d^i(t_d^{-1})E)\subset\varphi_d^{i-1}(t_d^{-1})\psi_d(E)\subset\varphi_d^{i-1}(t_d^{-1})E$$where
the second inclusion uses the assumption. Taking the product over all $d$ this implies\begin{gather}\psi_D(\varphi_D^i(t_D^{-1})E)\subset\varphi_D^{i-1}(t_D^{-1})\psi_D(E)\subset\varphi_D^{i-1}(t_D^{-1})E.\label{colmepeda}\end{gather}From $\varphi_d(t_d^{-1})=t_d^{-q}$ we get $$\psi_d(t_d^{-1}E)\subset\psi(\varphi_d(t_d^{-1})E)\subset t_d^{-1}E$$and taking the product over all $d$ we thus get$$\psi_D(t_D^{-1}E)\subset\psi(\varphi_D(t_D^{-1})E)\subset t_D^{-1}E.$$Moreover, if $n(x)\in{\mathbb N}$ is such that $x\in\varphi_D^n(t_D^{-1})E$ for $n\ge n(x)$,
then iterated application of formula (\ref{colmepeda}) shows $$\psi_D^n(x)\in
\psi_D^n(\varphi_D^n(t_D^{-1})E)\subset
\psi_D^{n-1}(\varphi_D^{n-1}(t_D^{-1})E)\subset\ldots\subset t_D^{-1}E$$for $n\ge n(x)$.\hfill$\Box$\\

 \begin{lem}\label{lotharjob} (a) There are lattices $E_0$, $E_1$ in ${\bf
    D}$ with $$\varphi_D(E_0)\subset t_DE_0\subset E_0\subset
   E_1\subset k[[t_{\bullet}]]\varphi_D(E_1).$$

(b) For
any $n\ge 0$ we have
$\psi_D^n(E_0)\subset\psi_D^{n+1}(E_0)\subset E_1$.
 \end{lem}

{\sc Proof:} (a) Let $e_1,\ldots, e_r$ be a generating system of ${\bf
    D}$ as a $k((t_{\bullet}))$-module. Then also $\varphi_d(e_1),\ldots,
\varphi_d(e_r)$ is such a generating system of ${\bf
    D}$, for each $d\in D$. We therefore find elements $f_{ij}$,
$g_{ij}$ in $k((t_{\bullet}))$ such that $\varphi_D(e_j)=\sum_{i=1}^rf_{ij}e_i$ and
$e_j=\sum_{i=1}^rg_{ij}\varphi_D(e_i)$ for all $1\le j\le r$. Choose
$n\in{\mathbb N}$ with $t_D^{n(q-1)}f_{ij}\in t_Dk[[t_{\bullet}]]$ and
$t_D^{n(q-1)}g_{ij}\in t_dk[[t_{\bullet}]]$ for all $i$, $j$. Then$$E_0=\sum_{i=1}^rt_D^nk[[t_{\bullet}]]e_i,\quad\quad\quad E_1=\sum_{i=1}^rt_D^{-n}k[[t_{\bullet}]]e_i$$work as
desired.

(b) Choose $a\in k[[t_{\bullet}]]$ with $\psi_D(a)=1$. For $x\in E_0$ we have
$\psi_D^n(x)=\psi_D^{n+1}(a\varphi_D(x))\in \psi_D^{n+1}(E_0)$ since
$\varphi_D(E_0)\subset t_D(E_0)$ implies $\varphi_D(x)\in E_0$ and hence
$a\varphi_D(x)\in E_0$. This shows
$\psi_D^n(E_0)\subset\psi_D^{n+1}(E_0)$. As $E_0\subset E_1\subset k[[t_{\bullet}]]\varphi_d(E_1)$, an induction using Lemma \ref{codoni} shows
$\psi_D^{n+1}(E_0)\subset E_1$.\hfill$\Box$\\
  
\begin{pro}\label{goldabreise} (a) There exists a unique lattice ${\bf
    D}^{\sharp}$ in ${\bf D}$ with $\psi_D({\bf D}^{\sharp})={\bf D}^{\sharp}$ and such that for each $x\in {\bf D}$ there is some $n\in{\mathbb N}$ with $\psi_D^n(x)\in {\bf D}^{\sharp}$.

(b) For any lattice $E$ in ${\bf D}$ we have $\psi_D^n(E)\subset {\bf D}^{\sharp}$ for all $n>>0$.

(c) For any lattice $E$ in ${\bf D}$ with $\psi_D(E)=E$ we have
$$t_D{\bf D}^{\sharp}\subset E\subset {\bf D}^{\sharp}.$$
\end{pro} 

{\sc Proof:} Let $E_0$, $E_1$ be as in Lemma \ref{lotharjob}. For
$n\in{\mathbb N}$ put $F_n=\psi_D^n(E_0)$. For $x\in E_0$ we have
$\psi_D^n(x)=\psi_D^{n+1}(\psi_D(1)\varphi_D(x))$; since $\psi_D(1)\varphi_D(x)\in E_0$ this
shows that $(F_n)_n$ is an increasing sequence of lattices in ${\bf D}$. As
$E_1\subset  k[[t_{\bullet}]]\varphi_D(E_1)$, Lemma \ref{codoni}
(c) shows $F_n\subset E_1$. As $k[[t_{\bullet}]]$ and hence $E_1$ is
noetherian, there is some $n_0$ with $F_n=F_{n_0}$ for all $n\ge n_0$, and
hence with $\psi_D(F_{n_0})=F_{n_0}$. For
$m\in{\mathbb N}$ put $$G_m=\psi_D^mt_D^{-1}F_{n_0}.$$Lemma \ref{codoni} (d) shows that 
$(G_m)_m$ is a descending sequence of lattices in ${\bf D}$, containing
$F_{n_0}$ since $\psi_D(F_{n_0})=F_{n_0}$. As $t_D^{-1}F_{n_0}/F_{n_0}$ is artinian we therefore find some $m_0$ with $G_m=G_{m_0}$ for all $m\ge m_0$, and
hence with $\psi_D(G_{m_0})=G_{m_0}$. Moreover, Lemma
\ref{codoni} (d) shows that for each $x\in G_{m_0}$ there is some
$i(x)\in{\mathbb N}$ with $\psi_D^i(x)\in t_D^{-1}F_{n_0}$ for all $i\ge
i(x)$. We then have $\psi_D^{m_0+i}(x)\in G_{m_0}$ for all $i\ge
i(x)$. Thus, ${\bf D}^{\sharp}=G_{m_0}$ works as desired.

To see the uniqueness of ${\bf D}^{\sharp}$, assume that there is another
candidate $\tilde{\bf D}^{\sharp}$ satisfying the same properties. Then so
does ${\bf D}^{\sharp}+\tilde{\bf D}^{\sharp}$, hence we may assume ${\bf
  D}^{\sharp}\subset \tilde{\bf D}^{\sharp}$. But $\psi_d$ for any $d\in D$ acts both
surjectively and nilpotently on the finite dimensional $k$-vector space $\tilde{\bf D}^{\sharp}/{\bf
  D}^{\sharp}$, hence ${\bf D}^{\sharp}=\tilde{\bf D}^{\sharp}$.\hfill$\Box$\\

\begin{pro}\label{wienabsage} (a) For any lattice $E$ in ${\bf D}$ contained
  in ${\bf D}^{\sharp}$ and stable under
  $\psi_d$ for $d\in D$ we have $\psi_d(E)=E$.

(b) The intersection ${\bf D}^{\natural}$ of all lattices in ${\bf D}$
  contained in ${\bf D}^{\sharp}$ and stable under
  $\psi_d$ for all $d\in D$ is itself a lattice, and it satisfies $\psi_d({\bf D}^{\natural})={\bf
    D}^{\natural}$ for all $d\in D$.
\end{pro} 

{\sc Proof:} (a) Since ${\bf D}^{\sharp}$ as well as $E$ and $\psi_d(E)$ are
lattices in ${\bf D}^{\sharp}$, both ${\bf D}^{\sharp}/E$ and ${\bf
  D}^{\sharp}/\psi_d(E)$ are finite dimensional $k$-vector spaces. $\psi_d$ induces an isomorphism ${\bf D}^{\sharp}/E={\bf D}^{\sharp}/\psi_d(E)$ (as $\psi_d(E)\subset E$), hence $\psi_d(E)=E$.

(b) For any lattice $E$ in ${\bf D}$
  contained in ${\bf D}^{\sharp}$ and stable under
  $\psi_d$ for all $d\in D$ we have $t_D{\bf D}^{\sharp}\subset E$ by what we saw
in (a) together with proposition
\ref{goldabreise}. This shows $t_D{\bf D}^{\sharp}\subset{\bf D}^{\natural}$, hence
${\bf D}^{\natural}$ is indeed a lattice, and $\psi_d({\bf D}^{\natural})={\bf
    D}^{\natural}$ follows by applying (a) once more.\hfill$\Box$\\

\begin{lem}\label{klara80} Let $\Delta$ be as in Lemma \ref{villpasteur}, with each
  $\varphi_{d}$ acting injectively on $\Delta$. If $\Delta[t_{\bullet}]$ generates $\Delta$ as a
  $k[[t_{\bullet}]][\varphi_{\bullet}]$-module then $\Delta^*={\bf D}^{\natural}$.  
\end{lem} 

{\sc Proof:} For $i=(i_d)_{d\in
  D}\in{\mathbb Z}_{\ge0}^D$ let $$F^i\Delta^*=\{\ell\in\Delta^*\,|\,\ell((\prod_{d\in
  D}t_{d}^{n_d}\varphi_{d}^{i_d})(x))=0\mbox{ for all
}n_d>0, x\in\Delta[t_{\bullet}]\}.$$This is a
$k[[t_{\bullet}]]$-submodule of $\Delta^*$. Let $E$ be a lattice in ${\bf D}$
contained in ${\bf D}^{\sharp}$ with
$\psi_{d}(E)\subset E$ for all $d\in D$. We have $\cap_iF^i\Delta^*=0$ since $\Delta[t_{\bullet}]$ generates
$\Delta$ as a
  $k[[t_{\bullet}]][\varphi_{\bullet}]$-module. As $E$
generates ${\bf D}$ we therefore find $F^i\Delta^*\subset E$ for some $i$. But $$(\prod_{d\in
  D}\psi_{d}^{i_d})F^i\Delta^*=\{\ell((\prod_{d\in D}\varphi_d^{i_d})(.))\,|\,\ell\in\Delta^*\}=\Delta^*$$where the second equality follows from the injectivity of the $\varphi_d$. We thus obtain $\Delta^*\subset E$ as the $\psi_{d}$ respect $E$.\hfill$\Box$\\

\begin{pro}\label{folklore} (Colmez) Suppose $|D|=1$. If ${\bf D}$ is an irreducible \'{e}tale $(\varphi,\Gamma)$-module with ${\rm dim}_{k((t))}({\bf D})\ge2$, then ${\bf D}^{\natural}={\bf D}^{\sharp}$. If ${\rm dim}_{k((t))}({\bf D})=1$ then ${\rm dim}_k({\bf D}^{\sharp}/{\bf D}^{\natural})=1$.
\end{pro}

{\sc Proof:} See \cite{col} Corollaire II.5.21 for the first statement. The second one follows e.g. from I.3.2. Exemple in \cite{colhaupt}, but also from the discussion of example (a) below.\hfill$\Box$\\

{\bf Remark:} It is easy to see that both ${\bf D}\mapsto {\bf D}^{\natural}$ and ${\bf D}\mapsto {\bf D}^{\sharp}$ are functors from the category of \'{e}tale $(\varphi_{\bullet},\Gamma_{\bullet})$-modules to the category of  $\psi_{\bullet}$-modules (obvious definition). Moreover, if \begin{gather}0\longrightarrow {\bf D}_1\longrightarrow {\bf D}\longrightarrow {\bf D}_2\longrightarrow0\label{heinrichbur50}\end{gather} is an exact sequence of \'{e}tale $(\varphi_{\bullet},\Gamma_{\bullet})$-modules, then the sequences \begin{gather}0\longrightarrow {\bf D}_1^{\natural}\longrightarrow {\bf D}^{\natural}\longrightarrow {\bf D}_2^{\natural}\longrightarrow0,\label{heinrich50}\\0\longrightarrow {\bf D}_1^{\sharp}\longrightarrow {\bf D}^{\sharp}\longrightarrow {\bf D}_2^{\sharp}\longrightarrow0\label{heinrichb50}\end{gather}are both exact on the left and on the right (see \cite{col} Proposition II 4.6 and  Proposition II 5.19 for the case $|D|=1$). However, in general they need not be exact in the middle. We are going to exemplify this below.

\section{Examples} 

{\bf (a)} By Proposition \ref{folklore}, if $|D|=1$ then a rank one \'{e}tale $(\varphi,\Gamma)$-module contains precisely two $(\psi,\Gamma)$-stable lattices with surjective $\psi$-operator. If $|D|>1$ we find more.

Fix some $c_{d}\in k^{\times}$ and some $m_{d}\in{\mathbb Z}/(q-1){\mathbb Z}$ for each $d\in D$. Put $$B=\bigoplus_{C\subset D}k.e_C,$$the $k$-vector space with basis $\{e_C\}_{C\subset D}$ indexed by the subsets $C$ of $D$. Let $k[[t_{\bullet}]][\Gamma_{\bullet}]$ act on $B$ by requiring\begin{align}t_{d}\cdot e_C=&\left\{\begin{array}{l@{\quad:\quad}l}0 &d\in C\\e_{C\cup\{d\}}&d\in D-C \end{array}\right.,\notag\end{align}\begin{align}\gamma_{d}\cdot e_C=&\left\{\begin{array}{l@{\quad:\quad}l}\gamma_{d}^{m_{d}+1}e_C &d\in C\\\gamma_{d}^{m_{d}}e_C&d\in D-C \end{array}\right..\notag\end{align}(On the right hand side of the defining formula for $\gamma_{d}\cdot e_C$ we refer to multiplication with the scalar in $k^{\times}$ to which $\gamma\in\Gamma={\mathcal O}_F^{\times}$ is mapped.)

Let ${\mathcal D}$ be a set of subsets of $D$ such that for any $C\in {\mathcal D}$ and $d\in D-C$ we also have $C\cup\{d\}\in {\mathcal D}$. It is clear that $\sum_{C\in {\mathcal D}}k.e_C$ is a $k$-sub vector space of $B$ stable under $k[[t_{\bullet}]][\Gamma_{\bullet}]$, hence $k[[t_{\bullet}]][\Gamma_{\bullet}]$ acts on $$B_{\mathcal D}=\frac{B}{\sum_{C\in {\mathcal D}}k.e_C}.$$Define $$\Delta_{\mathcal D}=\Delta_{\mathcal D}(c_{\bullet}, m_{\bullet})=\frac{k[[t_{\bullet}]][\varphi_{\bullet}]\otimes_{k[[t_{\bullet}]]}B_{\mathcal D}}{\langle t_{d}^{q-1}\varphi_{d}\otimes e_{\emptyset}-1\otimes c_{d}e_{\emptyset} \rangle_{d\in D}}$$where $\langle ? \rangle$ indicates the ${{k}}((t_{\bullet}))[\varphi_{\bullet}]$-sub module generated by all expressions within the brackets (and $e_{\emptyset}$ actually means the class of $e_{\emptyset}\in B$ in $B_{\mathcal D}$). One checks that this submodule $\langle t_{d}^{q-1}\varphi_{d}\otimes e_{\emptyset}-1\otimes c_{d}e_{\emptyset} \rangle_{d\in D}$ is in fact also stable under the action of $\Gamma_{\bullet}$; indeed, the $t_{d}^{q-1}\varphi_{d}\otimes e_{\emptyset}-1\otimes c_{d}e_{\emptyset}$ are eigenvectors under the action of $\Gamma_{\bullet}$. It follows that $\Delta_{\mathcal D}$ becomes a ${{k}}[[t_{\bullet}]][\varphi_{\bullet},\Gamma_{\bullet}]$-module. It is finitely
  generated over $k[[t_{\bullet}]][\varphi_{\bullet}]$, admissible over
  $k[[t_{\bullet}]]$ and each $\psi_d$ acts surjectively on $\Delta_{\mathcal D}^*=(\Delta_{\mathcal D})^*$. Thus $\Delta_{\mathcal D}^*$ is a lattice inside$${\bf D}={\bf
    D}(c_{\bullet}, m_{\bullet})=\Delta_{\mathcal D}^*\otimes_{k[[t{\bullet}]]}k((t_{\bullet})).$$The natural projections $B_{{\mathcal D}}\to B_{{\mathcal D}'}$ for ${\mathcal D}\subset{\mathcal D}'$ induce $k[[t_{\bullet}]]$-linear inclusions $\Delta_{{\mathcal D}'}^*\to  \Delta_{{\mathcal D}}^*$ which, when tensored with $k((t_{\bullet}))$, become isomorphisms. In particular, the rank one \'{e}tale
  $(\varphi_{\bullet},\Gamma_{\bullet})$-module ${\bf D}$ is in a natural way independent of ${\mathcal D}$, and (inside ${\bf
    D}$) we have$$\Delta_{\mathcal D}^*{\neq} \Delta_{{\mathcal D}'}^*\quad\quad\mbox{ whenever }\quad{\mathcal D}\neq{\mathcal D}'.$$Taking ${\mathcal D}_1$ to be set empty set, so that $B=B_{{\mathcal D}_1}$, we find $\Delta_{{\mathcal D}_1}^*={\bf D}^{\sharp}$. Taking ${\mathcal D}_0=\{D\subset C\,|\,D\ne \emptyset\}$, so that $B_{{\mathcal D}_0}$ is of $k$-dimension $1$ (generated by the class of $e_{\emptyset}$), we find $\Delta_{{\mathcal D}_0}^*={\bf D}^{\natural}$, cf. Lemma \ref{klara80}.\\

  In the following examples we choose (as we may) the coordinate $t$ such that $\Phi(t)=\pi t+t^q$. We assume $|D|=1$ and drop subscripts $(.)_{d}$. We describe various $\Delta$'s as quotients $\Delta=(k[[t]][\varphi]\otimes_kM)/\nabla$ with finite dimensional $k$-vector spaces $M$, and where always $\nabla$ is generated as a $k[[t]][\varphi]$-submodule by elements in $k[[t]]\varphi\otimes_kM+k[[t]]\otimes_kM$ only (i.e. no higher powers of $\varphi$ occur in these generators).

  In all these examples, $\varphi$ acts injectively in $\Delta$ (hence $\psi$ acts surjectively on $\Delta^*$) and $t$ acts surjectively on $\Delta$ (hence $\Delta^*$ is $t$-torsion free).\\
  
{\bf (b)} We describe a $\Delta$ defining an extension between two rank one \'{e}tale $(\varphi,\Gamma)$-modules.


Fix $\alpha\in k$. Let $\langle e_1, e_2,f\rangle_k$ denote the $k$-vector space with basis
$\{e_1, e_2,f\}$. In $k[[t]][\varphi]\otimes_k\langle e_1, e_2,f\rangle_k$ consider
the subset\footnote{In writing the elements of ${\mathcal R}$ we suppress the symbol $\otimes$.}\begin{align}{\mathcal R}=\{&te_1-e_2,\quad te_2,\quad tf,\notag\\{}&t^{q-1}\varphi f-f,\quad
t^{q-1}\varphi e_1-e_1-\alpha t^{q-2}\varphi f\}\notag\end{align}and let $\Delta$ denote the
quotient of $k[[t]][\varphi]\otimes_k\langle e_1, e_2,f\rangle_k$ by the
$k[[t]][\varphi]$-submodule generated by the elements in ${\mathcal R}$. Let $\langle f\rangle_k$ denote the $k$-sub vector space of $\langle e_1, e_2,f\rangle_k$ spanned by $f$ and define the $k$-vector space $\overline{\langle e_1, e_2\rangle_k}$ by the exact sequence$$0\longrightarrow \langle f\rangle_k\longrightarrow \langle e_1, e_2,f\rangle_k\longrightarrow\overline{\langle e_1, e_2\rangle_k}\longrightarrow0.$$We identify $\{e_1, e_2\}$ with a $k$-basis of $\overline{\langle e_1, e_2\rangle_k}$. We obtain an exact sequence$$0\longrightarrow
\frac{k[[t]][\varphi]\otimes_k\langle f\rangle_k}{\nabla_1}\longrightarrow
\Delta\longrightarrow\frac{k[[t]][\varphi]\otimes_k\overline{\langle e_1, e_2\rangle_k}}{\nabla_2}\longrightarrow
0$$where $\nabla_1$
(resp. $\nabla_2$) is the respective $k[[t]][\varphi]$-submodule generated by
$tf$ and $t^{q-1}\varphi f-f$ (resp. by $te_1-e_2, te_2$ and $t^{q-1}\varphi e_1-e_1$). 

Next, fix $a\in{\mathbb Z}$ and let $\gamma\in\Gamma$ act on $\langle e_1, e_2,f\rangle_k$ by means of $$\gamma\cdot
f=\gamma^{2+a} f\quad\quad\mbox{  and }\quad\quad\gamma\cdot e_i=\gamma^{i+a} e_i.$$(Here, in the expression $\gamma^{2+a} f$ resp. $\gamma^{i+a} e_i$ the $\gamma$ refers to the scalar in ${\mathbb F}_q^{\times}$ to which $\gamma\in\Gamma={\mathcal O}_F^{\times}$ is projected.) Then $$k[[t]][\varphi]\otimes_k\langle e_1, e_2,f\rangle_k \cong k[[t]][\varphi,\Gamma]\otimes_{k[\Gamma]}\langle e_1, e_2,f\rangle_k$$ so that $k[[t]][\varphi]\otimes_k\langle e_1, e_2,f\rangle_k$ receives a $k[[t]][\varphi,\Gamma]$-action. One checks that all
elements in ${\mathcal R}$ are eigenvectors for the action of $\Gamma$.  (For
the elements $t^{q-1}\varphi f-f$ and $t^{q-1}\varphi e_1-e_1-\alpha t^{q-2}\varphi f$ this computation
uses that $[\gamma](t)\equiv \gamma t$ modulo $t^qk[[t]]$, as is implied by our
assumption $\Phi(t)=\pi t+t^q$, see Lemma 0.1 in  \cite{elmarsusi}.) It follows
that $\Delta$, as well as the above exact sequence are in fact $k[[t]][\varphi,\Gamma]$-equivariant. 

In view of Proposition \ref{nopsi} we get an induced exact sequence (\ref{heinrichbur50}) of \'{e}tale
$(\varphi,\Gamma)$-modules, with ${\rm dim}_{k((t))}({\bf D}_1)= {\rm dim}_{k((t))}({\bf D}_2)=1$. If $F={\mathbb Q}_p$ then the Galois character attached (by Theorem \ref{sosego}) to ${\bf D}_1$ is obtained from the one attached to ${\bf D}_2$ by multiplying with the cyclotomic character. We have $\Delta={\bf D}^{\natural}$. Neither the sequence (\ref{heinrich50}) nor the sequence (\ref{heinrichb50}) is exact.\\

{\bf (c)} In contrast to what one might be tempted to think in view of Proposition \ref{folklore}, the possible failure of exactness of the sequences (\ref{heinrich50}) or (\ref{heinrichb50}) can {\it not} exclusively be reduced to the non-uniqueness of $(\psi,\Gamma)$-stable lattices with surjective $\psi$-operator inside \'{e}tale $(\varphi,\Gamma)$-module of rank one.

Let $\langle e_1, e_2, \tilde{e}, f_1, f_2\rangle_k$ denote the $k$-vector space with basis
$\{e_1, e_2, \tilde{e}, f_1, f_2\}$. Let $0\le s\le q-1$. In $k[[t]][\varphi]\otimes_{k}\langle e_1, e_2, \tilde{e}, f_1, f_2\rangle_k$ consider
the subset\begin{align}{\mathcal R}=\{&t^{q-1}\varphi
e_1-e_2-t^{s}\varphi f_2,\quad\varphi
e_2-e_1,\quad t^{q-2-s}\varphi
f_1-f_2,\quad t^{1+s}\varphi
f_2-f_1,\notag\\{}&t e_1-\tilde{e},\quad t\tilde{e},\quad t e_2,\quad t f_1,\quad t f_2\}\notag\end{align}and let $\Delta$ denote the
quotient of $k[[t]][\varphi]\otimes_{k}\langle e_1, e_2, \tilde{e}, f_1, f_2\rangle_k$ by the
$k[[t]][\varphi]$-submodule generated by the elements in ${\mathcal R}$. The exact sequence $0\to\langle f_1, f_2\rangle_k\to\langle e_1, e_2, \tilde{e}, f_1, f_2\rangle_k\to\langle e_1, e_2, \tilde{e}\rangle_k\to 0$ gives rise to an exact sequence$$0\longrightarrow
\frac{k[[t]][\varphi]\otimes_{k}\langle f_1, f_2\rangle_k}{\nabla_1}\longrightarrow \Delta\longrightarrow\frac{k[[t]][\varphi]\otimes_{k}\langle e_1, e_2, \tilde{e}\rangle_k}{\nabla_2}\longrightarrow 0,$$where $\nabla_1$
(resp. $\nabla_2$) is the respective $k[[t]][\varphi]$-submodule generated by
 $t^{q-2-s}\varphi
f_1-f_2$ and $t^{1+s}\varphi
f_2-f_1$ and $t f_1$, $t f_2$ (resp. by $t^{q-1}\varphi
e_1-e_2$ and $\varphi
e_2-e_1$ and $t e_1-\tilde{e}$, $t\tilde{e}$, $t e_2$). 

Next, fix $a\in{\mathbb Z}$ and let $\gamma\in\Gamma$ act on $\langle e_1, e_2, \tilde{e}, f_1, f_2\rangle_k$ by means of$$\gamma\cdot e_1=\gamma^ae_1,\quad\gamma\cdot e_2=\gamma^ae_2,\quad\gamma\cdot \tilde{e}=\gamma^{a+1}\tilde{e},\quad\quad\gamma\cdot f_1=\gamma^{1+a}f_1,\quad\gamma\cdot f_2=\gamma^{a-s}f_2.$$Then $$k[[t]][\varphi]\otimes_k \langle e_1, e_2, \tilde{e}, f_1, f_2\rangle_k\cong k[[t]][\varphi,\Gamma]\otimes_{k[\Gamma]}\langle e_1, e_2, \tilde{e}, f_1, f_2\rangle_k$$ so that $k[[t]][\varphi]\otimes_k\langle e_1, e_2, \tilde{e}, f_1, f_2\rangle_k$ receives a $k[[t]][\varphi,\Gamma]$-action. One checks that all
elements in ${\mathcal R}$ are eigenvectors for the action of $\Gamma$. It follows
that $\Delta$, as well as the above exact sequence are in fact $k[[t]][\varphi,\Gamma]$-equivariant. In view of Proposition \ref{nopsi} we get an induced exact sequence (\ref{heinrichbur50}) of \'{e}tale
$(\varphi,\Gamma)$-modules, with ${\rm dim}_{k((t))}({\bf D}_1)= {\rm dim}_{k((t))}({\bf D}_2)=2$. We have ${\bf D}_1^{\sharp}={\bf D}_1^{\natural}$
and ${\bf D}_2^{\sharp}={\bf D}_2^{\natural}$, with both ${\bf D}_1$ and ${\bf D}_2$ being irreducible, but the sequence (\ref{heinrichb50}) is not exact.\\

{\bf (d)} Let $\langle e_1, e_2, f_1, f_2\rangle_k$ denote the $k$-vector space with basis
$\{e_1, e_2, f_1, f_2\}$. We view it as a $k[[t]]$-module with trivial action by
$t$. Fix $0\le s\le k\le q-1$. In $k[[t]][\varphi]\otimes_{k[[t]]}\langle e_1, e_2, f_1, f_2\rangle_k$ consider
the subset\begin{align}{\mathcal R}=\{&t^{k}\varphi
e_1-e_2+t^{s}\varphi f_2,\quad t^{q-1-k}\varphi
e_2-e_1,\quad t^{k-s}\varphi
f_1-f_2,\quad t^{q-1-k+s}\varphi
f_2-f_1\}\notag\end{align}and let $\Delta$ denote the
quotient of $k[[t]][\varphi]\otimes_{k[[t]]}\langle e_1, e_2, f_1, f_2\rangle_k$ by the
$k[[t]][\varphi]$-submodule generated by the elements in ${\mathcal R}$. One
first checks that there is a natural exact sequence$$0\longrightarrow
\frac{k[[t]][\varphi]\otimes_{k[[t]]}\langle f_1, f_2\rangle_k}{\nabla_1}\longrightarrow \Delta\longrightarrow\frac{k[[t]][\varphi]\otimes_{k[[t]]}\langle e_1, e_2\rangle_k}{\nabla_2}\longrightarrow 0,$$where $\nabla_1$
(resp. $\nabla_2$) is the respective $k[[t]][\varphi]$-submodule generated by
 $t^{k-s}\varphi
f_1-f_2$ and $t^{q-1-k+s}\varphi
f_2-f_1$ (resp. by $t^{k}\varphi
e_1-e_2$ and $t^{q-1-k}\varphi
e_2-e_1$). 

Next, let $\Gamma$ act on $\langle e_1, e_2, f_1, f_2\rangle_k$ by means of
$$\gamma\cdot e_2=\gamma^ke_2,\quad \gamma\cdot f_2=\gamma^{k-s}f_2, \quad\gamma\cdot f_1=f_1,\quad \gamma\cdot e_1=e_1$$(understanding $\gamma^ke_2$ and $\gamma^{k-s}f_2$ similarly as before). Then $$k[[t]][\varphi]\otimes_{k[[t]]}\langle e_1, e_2, f_1, f_2\rangle_k\cong {k[[t]][\varphi,\Gamma]}\otimes_{k[[t]][\Gamma]}\langle e_1, e_2, f_1, f_2\rangle_k$$ so that $k[[t]][\varphi]\otimes_{k[[t]]}\langle e_1, e_2, f_1, f_2\rangle_k$ receives a $k[[t]][\varphi,\Gamma]$-action. Now one checks that all
elements in ${\mathcal R}$ are eigenvectors for the action of $\Gamma$. It follows
that $\Delta$, as well as the above exact sequence are in fact $k[[t]][\varphi,\Gamma]$-equivariant. In view of Proposition \ref{nopsi} we get an induced exact sequence (\ref{heinrichbur50}) of \'{e}tale
$(\varphi,\Gamma)$-modules, with ${\rm dim}_{k((t))}({\bf D}_1)= {\rm dim}_{k((t))}({\bf D}_2)=2$. It does not split. The  \'{e}tale
$(\varphi,\Gamma)$-module ${\bf D}$ lies in the essential image of the functor from supersingular modules over the pro-$p$ Iwahori Hecke algebra of ${\rm GL}_2(F)$ to \'{e}tale
$(\varphi,\Gamma)$-modules constructed in \cite{elmarsusi} if and only if $s=0$.

 \end{document}